\RequirePackage{fix-cm}
\documentclass[smallextended]{svjour3}       
\smartqed  
\usepackage{graphicx}
\usepackage{amsmath}
\usepackage{amssymb, comment}
\usepackage{mathrsfs}
\usepackage{multirow}
\usepackage[colorlinks=true]{hyperref}
\numberwithin{equation}{section}

\def \E {\mathbb E}

\def \P {\mathbb{P}}

\def \approxd {\,{\buildrel d \over \approxd}\,}

\newtheorem{prop}{Proposition}

\usepackage{subcaption}

\renewcommand{\proofname}{\noindent {\bf Proof}}

\title{A single server queue with batch arrivals and semi-Markov services}

\begin{document}

\author{Abhishek \and
        Marko~Boon \and
        Onno~Boxma \and
        Rudesindo~N\'u\~nez-Queija
        }

\institute{Abhishek $\&$  Rudesindo~N\'u\~nez~Queija \at
           Korteweg-de Vries Institute for Mathematics, University of Amsterdam, Amsterdam, The Netherlands\\
           \email{\{Abhishek, nunezqueija\}@uva.nl}
         \and
        Marko~Boon $\&$ Onno~Boxma \at
        Department of Mathematics and Computer Science, Eindhoven University of Technology, Eindhoven, The Netherlands \\
        \email{m.a.a.boon@tue.nl, o.j.boxma@tue.nl}
       }

%





\date{}
\maketitle
\begin{abstract}
We investigate the transient and stationary queue-length distributions of a class of service systems with correlated service times.
The classical $M^X/G/1$ queue with semi-Markov service times is the most prominent example in this class and serves as a vehicle to display our results.
The sequence of service times is governed by a modulating process $J(t)$.
The state of $J(\cdot)$ at a service initiation time determines the joint distribution of the subsequent service duration and the state of $J(\cdot)$ at the next service initiation.

Several earlier works have imposed technical conditions, on the zeros of a matrix determinant arising in the analysis, that are required in the computation of the stationary queue length probabilities.
The imposed conditions in several of these articles are difficult or impossible to verify.
Without such assumptions, we determine both the transient and the steady-state joint distribution of the number of customers immediately after a departure and the state of the process $J(t)$ at the start of the next service.

We numerically investigate how the mean queue length is affected by variability in the number of customers that arrive during a single service time.
Our main observations here are that increasing variability may {\em reduce} the mean queue length, and that the Markovian dependence of service times can lead to large queue lengths, even if the system is not in heavy traffic.

\keywords{batch arrivals, $M^X/G/1$ queue, semi-Markov service times, correlated service times, stationary and transient queue length analysis.}
\end{abstract}


\section{Introduction}
Service systems with correlated service durations have a long tradition in the queueing literature.
Such systems enjoy a large variety of application domains, including logistics, production management and telecommunications~\cite{AK,gaver,neuts66,takacs}.
Our main motivation stems from road traffic analysis, where traffic flows may interact at junctions or crossings~\cite{abhishekcomsnets2016,helbing2016}.
Focus, for illustration, on a traffic flow that merges into a main flow (very similar considerations are valid for road intersections).
If the traffic density on the main flow is high, vehicles in the secondary flow may queue up before merging into the main flow.
The merging times required for two subsequent vehicles will be strongly correlated as they experience similar traffic conditions on the main flow.
In this paper we will capture this dependence in a queueing model in which the sequence of service times is governed by a modulating Markovian process.
Although our analysis allows for a slightly larger class of models, we will use the classical $M/G/1$ queue with semi-Markov service times~\cite{neuts66}, and more specifically its extension to batch arrivals~\cite{neuts77a} to compare our results with existing literature.

The first to have investigated this class of queueing models was Gaver~\cite{gaver}, who derived the waiting time in a single-server queue with two types of customers arriving according to independent Poisson processes.
In that model, service times are class-specific and when service switches from one type to the other, an additional switch-over time is required.
This framework was generalized by Neuts~\cite{neuts66}, allowing for more than two customer types and the sequence of service times forming a semi-Markov process. Under technical assumptions (these will be discussed later in detail) Neuts obtained the transient and stationary distributions of queue lengths, waiting times and busy periods.
Subsequently, \c{C}inlar~\cite{cin_s} obtained the transient and stationary queue length distributions under less restrictive assumptions, and Purdue~\cite{pur} showed that the assumptions imposed by Neuts and \c{C}inlar are not necessary for the analysis of the busy period, presenting an alternative approach.
The literature on extensions of this model steadily expanded in the next two decades.
In~\cite{neuts77b}, Neuts studied the multi-type $M/G/1$ queue with change-over times when switching service from one type of customer to another. A further generalization allowing for Poisson arrivals of groups (batches) of customers of arbitrary random size was investigated by Neuts in~\cite{neuts77a}, obtaining the busy period, queue length and waiting time distributions. \\
The departure process of a related model with single Poisson arrivals and exponential service times was determined by Magalh\~aes and Disney~\cite{maga}.
In that model, the rate of the exponential service times depends on the type of the customer being served as well as that of its predecessor.\\
Models with single arrivals, but with both the arrivals and the services depending on a common semi-Markov process have been investigated by De Smit~\cite{desmit86} and Adan and Kulkarni~\cite{AK}.
Using the Wiener-Hopf factorization technique, De Smit~\cite{desmit86} obtained the waiting time and queue length distributions. Adan and Kulkarni~\cite{AK} considered a similar setting, but with the customer type being determined at arrival instants (independent of the service durations).\\

In this paper we investigate the transient and stationary queue length distributions in a single server model with semi-Markov service times and with batch arrivals (our framework includes Poisson arrivals of batches as the most prominent example).
In order to explain the technical contribution of our work, it is best to compare with the expositions of Neuts~\cite{neuts66} and \c{C}inlar~\cite{cin_s}.
In those papers only single Poisson arrivals were allowed, but the subsequent analysis is very similar.
The earlier mentioned technical assumptions made by Neuts entail that the zeros of a particular matrix determinant appearing in the transient analysis are either strictly separated or completely coincide. This ensures that the zeros are analytic functions of the entries of the matrix and, consequently, that the stationary distribution can be obtained from the transient distribution.
The assumptions were relaxed by \c{C}inlar~\cite{cin_s} while maintaining the analyticity of the zeros.
Unfortunately, it remains hard, if not impossible, to verify the required conditions in practice, as they must hold for the zeros as {\em functions} of the matrix entries.
As noted earlier, Purdue~\cite{pur} showed that the assumptions imposed by Neuts and \c{C}inlar are not necessary for the analysis of the busy period.
Our work show that these assumptions are not needed for the analysis of the queue length distribution either.
This comes at the expense of a separate analysis for the stationary distribution, which is more involved than that of the transient distribution.
{
Specifically, we determine the generating function of the number of customers immediately after the departure of an arbitrary customer, considering  both transient and steady-state behavior. For Poisson batch arrivals, in steady state we further obtain the queue length distribution at batch-arrival instants and at arbitrary times, which are identical due to PASTA. Note that this distribution is in general {\em not} the same as that at departure times (for single arrivals they would coincide).

A further contribution is an extensive numerical investigation of the mean queue length in steady state. We show that due to the dependence between service times, the mean number of customers may be very large, even if the load on the system is not large. A noteworthy observation is that {\em increasing} the variability in the number of customers arriving during a service time may in fact {\em decrease} the mean queue length. \\

The remainder of this paper is organized as follows. Section~\ref{sec:model} gives the model description in two layers. First we describe the $M^X/G/1$ queue with semi-Markov services and then present a somewhat more general framework. In Section~\ref{sec:queuedep}, we derive the transient and stationary probability generating functions of the number of customers in the system immediately after a departure. In Section~\ref{sec:queuear}, we derive the generating functions of the stationary number of customers at an arbitrary epoch, at batch arrival epochs and at customer arrivals.
The special case with only two customer types is specified in Section~\ref{sec:two}. Finally, in Section~\ref{sec:num}, we present numerical examples to demonstrate the impact of the correlated arrivals, and of the variability of the number of customers arriving during a service time, on the expected number of customers in the system.

\section{Model description}
\label{sec:model}

We start by describing the $M^X/G/1$ queueing model with semi-Markov services, which is the most natural example in our framework.
Our analysis extends directly to any model that satisfies the dynamics described in the recurrence relation~\eqref{recurA} below.

\subsection{The $M^X/G/1$ queue with semi-Markov service times}
Customers arrive in batches at a single server queue according to a Poisson process with rate $\lambda$;  the batch size is denoted by the random variable $B$ with generating function $B(z)$, for $|z|\leq1$. Customers are served in order of arrival, with speed $1$.
Customers within a batch are assumed to be ordered arbitrarily.
The service times are governed by a Markov process $J_n$, $n=0,1,\dots$, that can take values in $\{1,2,\dots,N\}$, for some integer $N$.
It will be convenient to refer to $J_n$ as the {\em type} of the $n$th customer; thus there are $N$ customer types.
The service time of the $n$th customer is denoted with $G^{(n)}$.
An essential feature of our model is that the type of the $(n+1)$th customer depends both on the type of the $n$th customer {\em and} on the service duration of the $n$th customer.
This exactly matches the framework of semi-Markov service times introduced by Neuts~\cite{neuts66}.
We define
\begin{equation}
G_{ij}(x) = \P(G^{(n)} \leq x,J_{n+1}=j|J_n=i), ~~~ x \geq 0, ~~~ i,j=1,2,\dots,N.
\label{Gij}
\end{equation}
For future use we introduce the Laplace-Stieltjes transform (LST)
\begin{equation}
\tilde{G}_{ij}(s) = \E[{\rm e}^{-sG^{(n)}} 1_{\{J_{n+1}=j\}}|J_n=i], ~~~ {\rm Re}~s \geq 0, ~~~ i,j=1,2,\dots,N,
\label{Gijs}
\end{equation} where $1_{\{.\}}$ denotes the indicator function.
In particular,
\begin{equation}
P_{ij} = G_{ij}(\infty) = \P(J_{n+1}=j|J_n=i), ~~~ i,j=1,2,\dots,N.
\label{Pij}
\end{equation}
The type of a customer, and its service time, do not depend on the arrival process.

It should be observed that $\{J_n, ~ n=1,2,\dots\}$ forms a finite-state Markov chain.
We shall restrict ourselves to irreducible Markov chains. The stationary distribution $\P(J=j)$ of the Markov chain $J_n$ is given by the unique solution of the set of equations
\begin{align}
\P(J=j) = \sum_{i=1}^N \P(J=i) P_{ij}, ~~~ j=1,2,\dots,N,
\label{Markov}
\end{align}
with normalizing condition $\sum_{j=1}^N \P(J=j) = 1$.
\\
The mean service time of an arbitrary customer is given by
\begin{equation}
\E [G] := \sum_{i=1}^N \sum_{j=1}^N \P(J_n=i) \E[G^{(n)} 1_{\{J_{n+1}=j\}}|J_n=i] .
\label{EB}
\end{equation}
The stability condition for this model is given by
\begin{equation}
\rho := \lambda \E[B]\E[G] < 1.
\label{stab}
\end{equation}
This can be formalized using Theorem 3 from Loynes~\cite{loynes62}, by describing the workload process in terms of ``super customers'' whose service times are the aggregate service times of customers in a single batch.
Let ${\mathcal G}^{(m)}$ be the service time of the super customer corresponding to the $m$th arriving batch, and ${\mathcal J}_{m}$ the type of the first customer in the $m$th batch.
Starting from a stationary version of the sequence $(G^{(n)},J_{n+1})$, one can readily construct a stationary sequence $({\mathcal G}^{(m)},{\mathcal J}_{m+1})$ for the super customers.
Note that by construction ${\mathcal G}^{(m)}$ is also stationary and, together with the arrival epochs of batches (which form an independent Poisson process), this sequence completely determines the workload process.
This description of the workload process satisfies the criteria to use the characterization for stability in Loynes~\cite{loynes62}.\\

\noindent
We will investigate the queue length process at departure times of customers. For that it will be convenient to define
$A_{n}$ as the number of customers arriving during the service time of the $n$th customer and $B_n$ as the size of the batch in which the $n$th customer arrived.
Note that for $i,j=1,2,\dots,N$, $|z| \leq 1$,
\begin{align}
A_{ij}(z):=\E[z^{A_{n}}1_{\{J_{n+1}=j\}}|J_n=i]= \tilde{G}_{ij}(\lambda(1-B(z))).
\label{AijzMXG1}
\end{align}
The queue length distribution at customer departure times is fully determined by the sequences $A_n$ and $B_n$.
For the analysis, it is not needed that the arrivals during service times occur in batches at Poisson instants.
For that reason we will now formulate our general model in terms of the $A_n$ and $B_n$ only; to specify our later results for the $M^X/G/1$ queue with semi-Markov services, we will simply substitute the relation given in~\eqref{AijzMXG1}.

\subsection{General model}
The inputs to our general model are probability generating functions of non-negative discrete random variables $A_{ij}(z)$, $i,j\in\{1,2,\dots,N\}$, and $B(z)$.
From the $A_{ij}(z)$, we construct a Markov process $(A_n,J_{n+1})$, $n=1,2,\dots$, satisfying
\begin{align}
\E[z^{A_{n}}1_{\{J_{n+1}=j\}}|J_n=i]= A_{ij}(z).
\label{GijzA}
\end{align}
In this construction it is implicit that $(A_n,J_{n+1})$ conditional on $J_n$ is independent of $A_{n-1}$.
The sequence $B_n$ is i.i.d.\ with generating function $B(z)$ and independent of the sequence $A_n$.

Next we define the recurrence relation
\begin{align}
X_n = \left\{
\begin{array}{l l}
\ X_{n-1}-1+A_n & \quad \text{if $X_{n-1} \geq 1$ }\\
A_n+B_n-1 & \quad \text{if $X_{n-1} =0$}
\end{array} \right., ~~~ n=1,2,3,\dots.
\label{recurA}
\end{align}

{\em Note:}
If the $A_{ij}(z)$ are set equal to \eqref{AijzMXG1}, then the sequence $X_n$ follows the same law as the number of customers at departure times in the $M^X/G/1$ queue with semi-Markov services.
The role of the $B_n$ is subtle in this representation: $B_n$ is only included if the $(n-1)$th customer leaves the system empty upon departure. The $n$th customer is therefore the first customer in a batch that arrives into an empty system.
Only for that reason, the sequence $B_n$ can be taken independent of the $A_n$ in the $M^X/G/1$ queue with semi-Markov services.\\

\par\noindent
In the sequel we will study the transient and stationary distributions of $X_n$ defined by \eqref{recurA}.
Again using Theorem 3 of Loynes~\cite{loynes62}, we may conclude that the stability condition in this case is
\begin{equation}
\rho:=\E[A]<1.
\label{eq:stabgen}
\end{equation}
Here $\E[A]$ denotes the expectation of the $A_n$ in stationarity:
\[
\E[A]=\sum_{i=1}^N\sum_{j=1}^N \P(J=i) \alpha_{ij},
\]
with
\begin{align}
\alpha_{ij}=\E[A_n 1_{\{J_{n+1}=j\}}|J_n=i]=A'_{ij}(1). \label{EA_ij}
\end{align}

Note that at first sight \eqref{recurA} does not seem to fit the framework in Loynes~\cite{loynes62}, because of the special condition when the system is empty. For stability, however, the behavior of an empty system is irrelevant.\\

\section{The queue length distribution at departure epochs}
\label{sec:queuedep}
We shall determine the transient and steady-state joint distribution of the number of customers immediately after a departure,
and the type of the next customer to be served.
From the recurrence relation \eqref{recurA} we find for the probability generating functions:
\begin{align}
&\E[z^{X_{n}}1_{\{J_{n+1}=j\}}]
=\E[z^{X_{n-1}-1+A_{n}}1_{\{J_{n+1}=j\}}1_{\{X_{n-1} \geq 1\}}]+\E[z^{A_n+B_n-1}1_{\{J_{n+1}=j\}}1_{\{X_{n-1}=0\}}]\nonumber\\
=&\E[z^{X_{n-1}-1+A_{n}}1_{\{J_{n+1}=j\}}]-\frac{1}{z}\E[z^{A_{n}}1_{\{J_{n+1}=j\}}1_{\{X_{n-1}=0\}}]+\E[z^{A_n+B_n-1}1_{\{J_{n+1}=j\}}1_{\{X_{n-1}=0\}}]\nonumber\\
=&\sum_{i=1}^{N}\E[z^{X_{n-1}-1+A_{n}}1_{\{J_{n+1}=j\}}|J_{n}=i]\P(J_{n}=i)-\frac{1}{z}\sum_{i=1}^{N}\E[z^{A_{n}}1_{\{J_{n+1}=j\}}1_{\{X_{n-1}=0\}}|J_{n}=i]\P(J_{n}=i)\nonumber\\
& +\sum_{i=1}^{N}\E[z^{A_n+B_n-1}1_{\{J_{n+1}=j\}}1_{\{X_{n-1}=0\}}|J_{n}=i]\P(J_{n}=i),\text{for} ~~ n=1,2,3,\dots, ~~~ j=1,2,\dots,N. \nonumber
\end{align}
Now we exploit the fact that $X_{n-1}$ and ($A_n,J_{n+1}$) are conditionally independent given $J_n$, and the $B_n$ are also independent of all other random variables:
\begin{align}
&\E[z^{X_{n}}1_{\{J_{n+1}=j\}}]
=\sum_{i=1}^{N}\E[z^{X_{n-1}-1}|J_{n}=i]\E[z^{A_{n}}1_{\{J_{n+1}=j\}}|J_{n}=i]\P(J_{n}=i)\nonumber\\
\ \ \ & +\frac{B(z)-1}{z}\sum_{i=1}^{N}\E[z^{A_{n}}1_{\{J_{n+1}=j\}}|J_{n}=i]\P(X_{n-1}=0|J_{n}=i)\P(J_{n}=i)\nonumber\\
=&\frac{1}{z}\sum_{i=1}^{N} \E[z^{X_{n-1}}1_{\{J_{n}=i\}}]\E[z^{A_{n}}1_{\{J_{n+1}=j\}}|J_{n}=i]\nonumber\\
\ \ \ & +\frac{B(z)-1}{z}\sum_{i=1}^{N}\E[z^{A_{n}}1_{\{J_{n+1}=j\}}|J_{n}=i]\P(X_{n-1}=0|J_{n}=i)\P(J_{n}=i),\nonumber\\
\ \ \ &  \text{for} ~~ n=1,2,3,\dots, ~~~ j=1,2,\dots,N. \nonumber\\
\label{twentysixAA}
\end{align}
\subsection{Steady-state analysis}
In this subsection, we restrict ourselves to the steady-state queue length distribution, assuming that the stability condition (\ref{eq:stabgen}) holds. In the next subsection, we will analyze the transient behavior of the queue length.\\
It will be useful to introduce some further notation:
for $i=1,2,\dots,N$,
\begin{align}
A_{i}(z)=\sum_{j=1}^{N}A_{ij}(z), \label{A_iz}
\end{align}
and,
\begin{align}
\alpha_i=\sum_{j=1}^{N}\alpha_{ij} \label{EA_i},
\end{align}
where the $\alpha_{ij}$ are defined in \eqref{EA_ij}.
Furthermore, for $j=1,2,\dots,N$, $|z| \leq 1$:
\begin{equation}
f_j(z)= {\rm lim}_{n \rightarrow \infty} \E[z^{X_n}1_{\{J_{n+1}=j\}}],
\label{fizA}
\end{equation}
\begin{equation}
f_j(0)= {\rm lim}_{n \rightarrow \infty} \P(X_n=0,J_{n+1}=j),
\label{fi0A}
\end{equation}
and note that,
\begin{equation}
f_j(1)= {\rm lim}_{n \rightarrow \infty} \P(J_{n+1}=j)=\P(J=j).
\label{fi1A}
\end{equation}
The probability generating function of the steady-state queue length distribution immediately after a departure is denoted by
\begin{equation}
F(z)=\sum_{j=1}^{N}f_j(z). \label{F(z)}
\end{equation}
In steady state, Equation \eqref{twentysixAA} leads to the following $N$ equations:
\begin{align}
& (z-A_{jj}(z))f_j(z)-\sum_{i=1,i\neq j}^{N}A_{ij}(z)f_i(z)
=(B(z)-1)\sum_{i=1}^{N}A_{ij}(z)
f_i(0),\quad \quad ~~~~~~ j=1,2,\dots,N. \label{twentysevenA}
\end{align}
We can also write these $N$ linear equations in matrix form as
\begin{align*}
M(z)^{T}f(z)=b(z),
\end{align*}
where \begin{align}
M(z)=&
\begin{bmatrix}
z-A_{11}(z) & -A_{12}(z) & \dots & -A_{1N}(z)\\
-A_{21}(z) & z-A_{22}(z) & \dots & -A_{2N}(z)\\
\dots & \dots & \dots & \dots \\
-A_{N1}(z) & -A_{N2}(z) & \dots & z-A_{NN}(z)
\end{bmatrix}, \label{Mmatrix}\\
f(z)=&
\begin{bmatrix}
f_1(z)\\
f_2(z)\\
\dots\\
f_N(z)
\end{bmatrix},
~~ b(z)=(B(z)-1)
\begin{bmatrix}
\sum_{i=1}^{N}A_{i1}(z)
f_i(0)\\
\sum_{i=1}^{N}A_{i2}(z)
f_i(0)\\
\dots \\
\sum_{i=1}^{N}A_{iN}(z)
f_i(0)
\end{bmatrix}. \nonumber
\end{align}
Therefore, solutions of the non-homogeneous linear system $M(z)^{T}f(z)=b(z)$ are in the form:
\begin{align}
f(z)=\frac{1}{\det M(z)^{T}}(\rm{cof}~M(z)^{T})^{T}b(z), ~~ \text{provided} ~ \det M(z) \neq 0. \label{M(z)}
\end{align}
Here $\rm{cof}~M(z)^T$ is the cofactor matrix of $M(z)^T$.
It remains to find the values of $f_1(0),f_2(0),\dots,f_N(0)$.
We shall derive $N$ linear equations for $f_1(0),f_2(0),\dots,f_N(0)$.
\\
{\bf First equation:}
\\
Note that $M(z)^{T}f(z)=b(z)$, which implies that
\begin{align*}
{\rm lim}_{z \rightarrow 1}\frac{1}{z-1}\hat{e}M(z)^{T}f(z)={\rm lim}_{z \rightarrow 1}\frac{1}{z-1}\hat{e}b(z),
 \end{align*} where $\hat{e}$ is a row vector with all entries one.\\
After simplification, we can write this as
\begin{align*}
{\rm lim}_{z \rightarrow 1}\frac{\sum_{i=1}^{N}\Big(z-\sum_{j=1}^{N}A_{ij}(z)\Big)f_i(z)}{z-1}={\rm lim}_{z \rightarrow 1}\frac{B(z)-1}{z-1}\sum_{j=1}^{N}\sum_{i=1}^{N}A_{ij}(z)f_i(0).
\end{align*}
Using $\sum_{i=1}^{N}f_i(1)=1$ and $\sum_{i=1}^{N}f_i(1)\alpha_i=\rho$, and after simplification, we get,
  \begin{align}
\sum_{i=1}^{N}f_i(0)=\frac{1-\rho}{\E[B]}. \label{thirtyfirstAA}
\end{align}
\\
{\bf (N-1) remaining equations:}
\\
To find the remaining $N-1$ equations, we first prove that $\det M(z)$ has exactly
$N-1$ zeros in $|z|<1$ and the zero $z=1$ on $|z|=1$. Since $f_i(z)$ is an analytic function in $|z|<1$, the numerator of $f_i(z)$ also has $N-1$ zeros in the unit disc $|z|<1$. As a consequence, these $N-1$ zeros provide $N-1$ linear equations for $f_1(0),f_2(0),\dots,f_N(0).$\\

To find the $N-1$ zeros, we use a method that has also been applied in \cite{AK,GHT,DeSmit}.
It is based on the concept of (strict) diagonal dominance in a matrix.
The proof consists of $4$ steps:
\\
Step 1: Prove that each element on the diagonal of $M(z)$ has exactly one zero in $|z|<1$.
\\
Step 2: Introduce a matrix $M(t,z)$, $0\leq t \leq 1$, with $M(1,z) = M(z)$,
and prove strict diagonal dominance of $M(t,z)$, i.e., each diagonal element of $M(t,z)$ is in absolute value larger
than the sum of the absolute values of the non-diagonal terms in the same row of the matrix.
\\
Step 3: Prove that $\det M(t,z)$ has exactly $N$ zeros in $|z|<1$ and none on $|z|=1$ for $0 \leq t < 1$.
\\
Step 4: Use continuity of $\det M(t,z)$ in $t$ for $0 \leq t < 1$ to prove that, indeed, $\det M(z)$ has $N-1$ zeros
in $|z|<1$ and one zero $z=1$ on $|z|=1$.
\\

\noindent
{\bf Step 1:} Prove that each element on the diagonal of $M(z)$ has exactly one zero in $|z|<1$.
\\
It follows from \eqref{Mmatrix} that $M(z)=D(z)+O(z)$,
where $D(z)$ is the diagonal matrix
\begin{align}
D(z)=
\begin{bmatrix}
z-A_{11}(z) & 0& \dots &0\\
0 & z-A_{22}(z) & \dots &0\\
\dots & \dots & \dots & \dots \\
0&0 & \dots & z-A_{NN}(z)
\end{bmatrix}
,
\label{Dmatrix}
\end{align}
and $O(z)$ is the off-diagonal matrix which corresponds to $M(z)$.
\\
\begin{prop}
$\det D(z)$ has exactly $N$ zeros (counting multiplicities) in $|z|<1$ and none satisfying $|z|=1$.
\label{Prop1}
\end{prop}
{\em Proof.}
First observe that
$\det D(z)= \prod_{i=1}^N (z-A_{ii}(z))$.
Because $|\frac{A_{ii}(z)}{z}|\leq P_{ii}<1$ on $|z|=1$,
Rouch\'e's theorem implies that the numbers of zeros of $z$ and $z-A_{ii}(z)$ are the same in $|z|<1$. $z$ has exactly one zero in $|z|<1$, and hence $z-A_{ii}(z)$ also has exactly one zero in $|z|<1$, for $i=1,2,\dots,N$.
\\
On $|z|=1$,
$|z-A_{ii}(z)|$ has no zeros, because
$|z-A_{ii}(z)| \geq |z| - |A_{ii}(z)| \geq 1-P_{11}>0$.
\\
Hence $\det D(z)$  has $N$ zeros in $|z|<1$ and none on $|z|=1$.
\\

Now we define the matrix
$M(t,z):=D(z)+tO(z)$,
where $0\leq t\leq 1$ is a real parameter.
Note that $M(0,z)=D(z)$ and $M(1,z)=M(z)$.\\

\noindent
{\bf Step 2:} Prove diagonal dominance for matrix $M(t,z)$.
\begin{prop}
$\det M(t,z) \neq 0$ for $0\leq t<1, |z|=1$ and for $t=1, |z|=1, z \neq 1$.
\label{Prop2}
\end{prop}
{\em Proof.}
Consider an arbitrary $i \in \{1,2,\dots,N\}$.
\begin{align}
|z-A_{ii}(z)|&\geq |z|-|A_{ii}(z)|
\nonumber
\\
&\geq 1-P_{ii}= \sum_{j \neq i} P_{ij}
>t \sum_{j \neq i} P_{ij}
\ \ \ \ \ \ \ {\rm for} ~  0\leq t<1, |z|=1.
\label{ineq}
\end{align}
On the other hand, $\sum_{j \neq i} |tA_{ij}(z)| \leq t \sum_{j \neq i} P_{ij}$
for $0\leq t<1, |z|=1$.
\\
Therefore, $|z-A_{ii}(z)|>|t \sum_{j \neq i} A_{ij}(z)|$ for  $0\leq t<1, |z|=1$.
This holds for $i=1,2,\dots,N$.
\\
Thus, $M(t,z)$ is strictly diagonally dominant. This implies that  $M(t,z)$ is a non-singular matrix,
i.e., $\det M(t,z) \neq 0$,
for $0\leq t<1, |z|=1$.
This concludes the proof for the case $0 \leq t < 1$, with $|z|=1$.
\\
We next turn to the case $t=1$, $|z|=1$, $z \neq 1$, again
considering an arbitrary $i \in \{1,2,\dots,N\}$.
Now (\ref{ineq}) is replaced by
$|z-A_{ii}(z)| >  \sum_{j \neq i} P_{ij}$ for  $ |z|=1, z\neq 1$.
On the other hand, $\sum_{j \neq i} |A_{ij}(z)| < \sum_{j \neq i} P_{ij}$.
Therefore, $|z-A_{ii}(z)|>| \sum_{j \neq i} A_{ij}(z)|$ for $|z|=1$, $z \neq 1$.
This holds for $i=1,2,\dots,N$.
In this way we have proven the strict diagonal dominance, and hence the non-singularity, also for $t=1$, $|z|=1$, $z \neq 1$.
\\

\noindent
{\bf Step 3:} Prove that $\det M(t,z)$ has exactly $N$ zeros in $|z|<1$ and none on $|z|=1$ for $0 \leq t < 1$.
\begin{prop}
\label{Prop3}
The function $\det M(t,z)$ has exactly $N$ zeros in $|z|<1$ and none on $|z|=1$ for $0\leq t<1$.
\end{prop}
{\em Proof.} Let $n(t)$ be the number of zeros of $\det M(t,z)$ in $|z|<1$.
By the argument principle, see Evgrafov~\cite[p.~97]{evgrafov},
\begin{equation}
n(t)=\frac{1}{2\pi i}\int_{|z|=1}\frac{\frac{\partial}{\partial z} \det M(t,z)}{\det M(t,z)} dz,
\end{equation}
where it should be noticed that $\det M(t,z)\neq 0$ on $|z|=1$ for $0\leq t<1$ according to Proposition~\ref{Prop2}.
Here, $n(t)$ is a continuous integer-valued function of $t$ for $0\leq t<1$ and $n(0)=N$ according to Proposition~\ref{Prop1}. So $n(t)=n(0)=N.$\\

From the above we may conclude that
$\det M(1,z)=M(z)$ has {\em at least} $N$ zeros in the closed unit disc,
because the zeros of $\det M(t,z)$ are continuous functions for $0\leq t\leq1$.
Finally we need to prove that there are {\em exactly} $N$ zeros in $|z| \leq 1$, one of which ($z=1$) lies on $|z|=1$.\\

\noindent
{\bf Step 4:} Use continuity of $\det M(t,z)$ in $t$ for $0 \leq t \leq 1$ to prove that $\det M(z)$ has $N-1$ zeros
in $|z| < 1$ and one zero $z=1$ on $|z|=1$.
\begin{prop}
$\frac{d}{dz}\{det M(z)\}|_{z=1}>0$ and $z=1$ is a simple zero of $\det M(z)$.
\label{Prop4}
\end{prop}
{\em Proof.}
Firstly, $z=1$ is a zero of $\det M(z)$.
Now we show that it is a simple zero.
Use that $\lim_{z\to 1} \frac{det ~ M(z)}{z-1}=\frac{d}{dz}\{det ~ M(z)\}|_{z=1}>0$, where the inequality is a consequence of the stability condition.
Hence, $z=1$ is a simple zero of  $\det M(z)$.
\begin{prop}
det  $M(t,1)>0$ for $0\leq t<1$.
\label{Prop5}
\end{prop}
{\em Proof.}
We shall exploit the fact that $\det M(t,1)$ is the product of all eigenvalues of $M(t,1)$.
So we need to prove that the product of these eigenvalues is positive.
\\
Consider the matrix $I - M(t,1)$, where $I$ is the identity matrix:
\\
\begin{align*}
I-M(t,1)=&
\begin{bmatrix}
P_{11} &t P_{12}&t P_{13}& \cdots& tP_{1N}\\
t P_{21} & P_{22} &t P_{23}&\cdots& tP_{2N}\\
t P_{31}&t P_{32} & P_{33}&\cdots& tP_{3N}\\
\vdots&\vdots&\vdots&&\vdots\\
tP_{N1}&tP_{N2}&tP_{N3}&&P_{NN}
\end{bmatrix}.
\end{align*}
Note that $I-M(t,1)$ is a substochastic matrix, so every eigenvalue of the matrix $I-M(t,1)$ lies in $|z|< 1$. Hence every eigenvalue of the matrix $M(t,1)$ lies in $|z-1|< 1$.
$M(t,1)$ is a real matrix, so if $M(t,1)$ has a complex eigenvalue, then the conjugate of this complex eigenvalue is also one of the eigenvalues of $M(t,1)$.
This implies that if $M(t,1)$  has complex eigenvalues, then the product of these complex eigenvalues is positive. The product of the real eigenvalues is also positive because every eigenvalue of the matrix $M(t,1)$ lies in $|z-1|< 1.$
This concludes the proof.
\\
\begin{prop}
The function  $\det M(z)$ has exactly $N-1$ zeros in $|z|<1$ and one zero on $|z|=1$ $(\text{at}  \ z=1)$.
\label{Prop6}
\end{prop}
{\em Proof.} We follow the argument of Gail et al.~\cite[p.~372]{GHT}. By letting $t \to 1$ in Proposition~\ref{Prop3}, it follows that $\det M(z)$ has at least $N$ zeros in $|z|\leq 1$.
By Proposition~\ref{Prop4}, given $\epsilon >0$, there is a real $z'$, $1-\epsilon <z'<1$, such that $\det M(z')$ is negative.
By continuity, there is a real $t'$, $1-\epsilon <t'<1$, such that $\det M(t',z')$ is negative.
Since $\det M(t',1)$ is positive according to Proposition~\ref{Prop5}, there is a real $z''$, $z'<z''<1$ with $\det M(t',z'')=0$.
Thus, the zero of $\det M(z)$ at $z=1$ is the limit of a zero of $\det M(t,z)$ from inside the unit disc.
As $t\to 1$, the limiting positions of the $N$ zeros of $\det M(t,z)$ are:  one at $z=1$ and the other $N-1$ in $|z|<1$.\\

\subsection{Transient analysis}
In this subsection, we shall determine the transient behavior of the probability generating function of the number of customers. The analysis proceeds largely analogous to the stationary case.
In fact, for the transient analysis, it turns out to be less involved to demonstrate the location of the roots.
We define
\begin{equation}
f_j(r,z)=\sum_{n=0}^{\infty}r^n \E[z^{X_{n}}1_{\{J_{n+1}=j\}}] ~~~~ \text{for}~~~~ |r|<1, j=1,2,...,N,\label{fjrzA}
\end{equation}
so that,
\begin{equation}
f_j(r,0)=\sum_{n=0}^{\infty}r^n\P(X_n=0,J_{n+1}=j).\label{fjrzeroA}
\end{equation}
Using \eqref{twentysixAA} with $\E[z^{A_{n}}1_{\{J_{n+1}=j\}}|J_n=i]=A_{ij}(z)$ in \eqref{fjrzA}, we get
\begin{align*}
f_j(r,z)=&\E[z^{X_{0}}1_{\{J_{1}=j\}}]+\frac{1}{z}\sum_{i=1}^{N}A_{ij}(z)\sum_{n=1}^{\infty}r^n \E[z^{X_{n-1}}1_{\{J_{n}=i\}}]\nonumber\\
\ \ \ & +\left(\frac{B(z)-1}{z}\right)\sum_{i=1}^{N}A_{ij}(z)\sum_{n=1}^{\infty}r^n\P(X_{n-1}=0,J_{n}=i)\\
=&z^{x_{0}}\P(J_{1}=j)+\frac{1}{z}\sum_{i=1}^{N}A_{ij}(z)\sum_{n=0}^{\infty}r^{n+1} \E[z^{X_{n}}1_{\{J_{n+1}=i\}}]\nonumber\\
\ \ \ & +r\left(\frac{B(z)-1}{z}\right)\sum_{i=1}^{N}A_{ij}(z)f_i(r,0),
\end{align*}
provided the initial number of customers in the system is deterministic and equal to $x_{0}$.\\

Using \eqref{fjrzA} and after simplification, we get the following $N$ equations:
\begin{align}
 (z-rA_{jj}(z))f_j(r,z)-r\sum_{i=1,i\neq j}^{N}A_{ij}(z)f_i(r,z)
=z^{X_{0}+1}\P(J_{1}=j)\nonumber\\
 +r\left(B(z)-1\right)\sum_{i=1}^{N}A_{ij}(z)f_i(r,0),\quad \quad ~~~~~~ j=1,2,\dots,N. \label{twentysevenAA}
\end{align}
We can also write these $N$ linear equations in matrix form as
\begin{align*}
M(r,z)^{T}f(r,z)=b(r,z),
\end{align*}
where
\begin{align*}
M(r,z)=&
\begin{bmatrix}
z-rA_{11}(z) & -rA_{12}(z) & \dots & -rA_{1N}(z)\\
-rA_{21}(z) & z-rA_{22}(z) & \dots & -rA_{2N}(z)\\
\dots & \dots & \dots & \dots \\
-rA_{N1}(z) & -rA_{N2}(z) & \dots & z-rA_{NN}(z)
\end{bmatrix},\\
f(r,z)=&
\begin{bmatrix}
f_1(r,z)\\
f_2(r,z)\\
\dots\\
f_N(r,z)
\end{bmatrix},~~
 b(r,z)=z^{X_{0}+1}\begin{bmatrix}
 \P(J_{1}=1)\\
 \P(J_{1}=2)\\
 \dots\\
 \P(J_{1}=N)
 \end{bmatrix}
 +r(B(z)-1)
\begin{bmatrix}
\sum_{i=1}^{N}A_{i1}(z)
f_i(r,0)\\
\sum_{i=1}^{N}A_{i2}(z)
f_i(r,0)\\
\dots \\
\sum_{i=1}^{N}A_{iN}(z)
f_i(r,0)
\end{bmatrix}.
\end{align*}
Therefore, solutions of the non-homogeneous linear system $M(r,z)^{T}f(r,z)=b(r,z)$ are in the form:
\begin{align}
f(r,z)=\frac{1}{\det M(r,z)^{T}}(\rm{cof}~M(r,z)^{T})^{T}b(r,z), ~~ \text{provided} ~ \rm{det} ~ M(r,z) \neq 0. \label{M(r,z)}
\end{align}
It remains to find the values of $f_1(r,0),f_2(r,0),\dots,f_N(r,0)$.
We shall derive $N$ linear equations for $f_1(r,0),f_2(r,0),\dots,f_N(r,0)$.\\

To find $N$ linear equations for $f_1(r,0),f_2(r,0),\dots,f_N(r,0)$, we fist prove that $\det M(r,z)$ has exactly $N$ zeros for fixed $r$ in $|z|<1$.
Since $M(r,z)=zI-rA(z)$, $\det M(r,z)$ is a continuous function in $r$ for $0\leq r \leq1$, and therefore the zeros are continuous in $0\leq r \leq1$.

{\em Remark.}
It is worth emphasizing that it is at this point that our approach is different from the analysis by Neuts~\cite{neuts66} and \c{C}inlar~\cite{cin_s}.
We do not require for each pair of elementary roots that they either be strictly different for all values of $0\leq r\leq1$ or coincide for all $0\leq r\leq1$.
The main price to pay is that we can not use that the roots are analytic in $r$ and we can therefore not obtain the stationary distribution from the transient distribution as $r\to1$.\\

\par\noindent

Compared to the steady-state analysis, the proof is simpler and only consists of two steps: \\

\noindent
{\bf Step 1:} Prove diagonal dominance of the matrix $M(r,z)$.
\begin{prop}
$\det M(r,z) \neq 0$ for $0\leq r <1, |z|=1$.
\label{Propr1}
\end{prop}
{\em Proof.}
Consider an arbitrary $i \in \{1,2,\dots,N\}$.
\begin{align}
|z-rA_{ii}(z)|&\geq |z|-r|A_{ii}(z)|
\nonumber
\\
&> 1-P_{ii}= \sum_{j \neq i} P_{ij}
> r \sum_{j \neq i} P_{ij}
\ \ \ \ \ \ \ {\rm for} ~  0\leq r<1, |z|=1.
\label{ineqr}
\end{align}
On the other hand, $\sum_{j \neq i} |rA_{ij}(z)| \leq r \sum_{j \neq i} P_{ij}$
for $0\leq r<1, |z|=1$.
\\
Therefore, $|z-r A_{ii}(z)|>|r \sum_{j \neq i} A_{ij}(z)|$ for  $0\leq r<1, |z|=1$.
This holds for $i=1,2,\dots,N$.
\\
Thus, $M(r,z)$ is strictly diagonally dominant. This implies that  $M(r,z)$ is a non-singular matrix,
i.e., $\det M(r,z) \neq 0$,
for $0\leq r<1, |z|=1$.
This completes the proof.\\

\noindent
{\bf Step 2:} Prove that $\det M(r,z)$ has exactly $N$ zeros in $|z|<1$ for $0 \leq r < 1$.
\begin{prop}
\label{Propr2}
The function $\det M(r,z)$ has exactly $N$ zeros in $|z|<1$ for $0\leq r<1$.
\end{prop}
{\em Proof.} Let $n(r)$ be the number of zeros of $\det M(r,z)$ in $|z|<1$.
As before, by the argument principle~\cite[~p.~97]{evgrafov},
\begin{equation}
n(r)=\frac{1}{2\pi i}\int_{|z|=1}\frac{\frac{\partial}{\partial z} \det M(r,z)}{\det M(r,z)} dz,
\end{equation}
where it should be noticed that $\det M(r,z)\neq 0$ on $|z|=1$ for $0\leq r<1$ according to Proposition~\ref{Propr1}.
Here, $n(r)$ is a continuous integer-valued function of $r$ for $0\leq r<1$ and $n(0)=N$ because $\det M(0,z)=z^n$ . So $n(r)=n(0)=N.$

\section{Poisson batch arrivals: stationary queue length at arrival and arbitrary epochs}
\label{sec:queuear}
In the previous section, we determined the stationary and the transient queue length distributions at departure times of customers. In the general framework, the exact arrival process of customers is not specified, but for the model with Poisson batch arrivals, we can obtain the stationary queue length distribution at {\em arbitrary time}, at {\em batch arrival instants} and at {\em customer arrival instants}.
Because of PASTA, the distribution of the number of customers already in system just before a new batch arrives (let us denote this by a generic random variable $X^{ba}$) coincides with the distribution of the number of customers in the system at an arbitrary time ($X^{arb}$).
The number of customers at customer arrival instants (denoted with $X^{ca}$) needs to be further specified, because with batch arrivals all customers in the same batch have the same arrival time.
As noted previously, customers within one batch are assumed to be (randomly) ordered.
Although they arrive at the same time, they see different numbers of customers in front of them.
In particular, the last customer in a batch sees all the customers that were already in the system {\em plus} all  other customers (excluding him/her) arriving in the same batch.
In the {\em customer average} distribution at arrival times, this must be taken into account.
In Figure~\ref{fig:levelcrossing} we depict three batch arrivals, two of which contain multiple customers and thus coincide with more than one customer arrival.
\begin{figure}[!ht]
\centering
\includegraphics[width=0.6\linewidth]{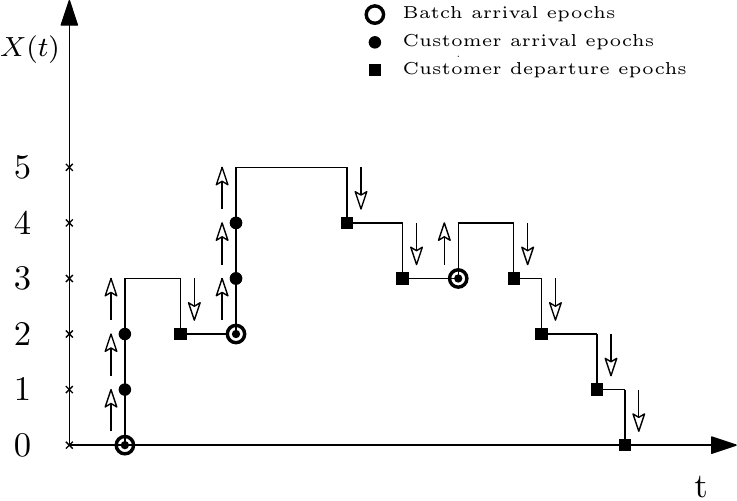}
\caption{Up and down-crossing.}
\label{fig:levelcrossing}
\end{figure}
Applying a simple level crossing argument with the aid of Figure~\ref{fig:levelcrossing}, it is readily seen that the distributions of $X$ (at departure times) and $X^{ca}$ must coincide: indeed, for each level $k=1,2,\dots$, customer departures that decrease the queue length from $k$ to $k-1$ must be matched by customer arrivals increasing the level from $k-1$ to $k$ (since the arrival of each customer within a batch is counted separately, the difference can be at most 1 which is negligible in the long run).\\
We can also link the distributions of $X^{ba}$ and $X^{ca}$:
A customer in an arriving batch sees in front of him the number of customers already in system ($X^{ba}$) and the number of customers in front of him in the same batch.
For an arbitrary customer in the batch the number of customers in front of him in the same batch has the forward recurrence distribution of $B$.
Summarizing:
 \begin{align}
 \E[z^{X}]=\E[z^{X^{ca}}] = \E[z^{X^{ba}}]\frac{1-B(z)}{\E[B](1-z)},
  \end {align}
  where we use independence of the batch size and the number of customers already in system,
and
 \begin{align}
 \E[z^{X^{arb}}]=\E[z^{X^{ba}}] .
 \label{F(z):arbitrary}
  \end {align}
From these relations we can obtain all the required distributions. It can be verified that these distributions agree with the results from Chaudhry\cite{chaudhary} for the model without dependencies between successive service times.

\section{The queueing model with  two customer types : departure epochs}
\label{sec:two}
In this section, we restrict ourselves to the case of two customer types, i.e., $N=2$. In this case, we are able to give an explicit expression for  the probability generating function of the number of customers in the system immediately after a departure.
For the steady-state behavior it follows from \eqref{twentysevenA} that:
\begin{align}
f_1(z)=\frac{\Big(B(z)-1\Big)\Big(f_1(0)\left(zA_{11}(z)+ A_{12}(z)A_{21}(z)-A_{11}(z)A_{22}(z) \right)+zf_2(0)A_{21}(z)\Big)}{(z-A_{11}(z))(z-A_{22}(z))-A_{12}(z)A_{21}(z)}, \label{twentynine}
\end{align}
\begin{align}
f_2(z)=\frac{\Big(B(z)-1\Big)\Big(zf_1(0) A_{12}(z)+f_2(0) \left(zA_{22}(z)+A_{12}(z)A_{21}(z)-A_{11}(z)A_{22}(z)\right) \Big)}{(z-A_{11}(z))(z-A_{22}(z))-A_{12}(z)A_{21}(z)}, \label{thirty}
\end{align}
where
\begin{equation}
f_1(0) = \frac{1-\rho}{\E[B]} \frac{A_{11}(\hat{z}) - \hat{z}}{A_{11}(\hat{z}) + A_{12}(\hat{z})-\hat{z}}, ~~~
f_2(0) = \frac{1-\rho}{\E[B]} \frac{A_{22}(\hat{z}) - \hat{z}}{A_{21}(\hat{z}) + A_{22}(\hat{z})-\hat{z}},\label{f1f2}
\end{equation} such that $f_1(0)+f_2(0)=\frac{1-\rho}{\E[B]}$ and $z=\hat{z}$ is the zero of $(z-A_{11}(z))(z-A_{22}(z))-A_{12}(z)A_{21}(z)$ with $|\hat{z}|<1$.\\

It is noted that the probability generating function of $X_n$ in steady state is
\begin{equation*}
F(z)= {\rm lim}_{n \rightarrow \infty} \E[z^{X_n}] .
\end{equation*}
 From Equation \eqref{F(z)}, for $N=2$, we can write $F(z)$ as the sum of $f_1(z)$ and $f_2(z)$, i.e.,
\begin{equation*}
F(z)=f_1(z)+f_2(z).
\end{equation*}
After substituting the values of $f_1(z)$ and $f_2(z)$ from Equations \eqref{twentynine} and \eqref{thirty} respectively, we obtain $F(z)$ as

\begin{align*}
F(z)=\frac{z(B(z)-1)\Big(f_1(0)( A_{11}(z)+A_{12}(z))+f_2(0) (A_{21}(z)+A_{22}(z))\Big)}{(z-A_{11}(z))(z-A_{22}(z))-A_{12}(z)A_{21}(z)} \nonumber \\
+ \frac{(B(z)-1)(f_1(0)+f_2(0))\Big(A_{12}(z)A_{21}(z)-A_{11}(z)A_{22}(z) \Big)}{(z-A_{11}(z))(z-A_{22}(z))-A_{12}(z)A_{21}(z)}.
\end{align*}

Equation \eqref{A_iz} states that $A_{i}(z)=A_{i1}(z)+A_{i2}(z)$ for $i=1,2.$
After substituting the values of $f_{i}(0)$ and $A_{i}(z)$ for $i=1,2$, $F(z)$ becomes
\begin{align*}
F(z)=\frac{z(B(z)-1)(1-\rho)\Big(c_1 A_{1}(z)+c_2A_{2}(z)\Big)}{\E[B]\Big((z-A_{11}(z))(z-A_{22}(z))-A_{12}(z)A_{21}(z)\Big)} \nonumber \\
+ \frac{(B(z)-1)(1-\rho)\Big(A_{12}(z)A_{21}(z)-A_{11}(z)A_{22}(z) \Big)}{\E[B]\Big((z-A_{11}(z))(z-A_{22}(z))-A_{12}(z)A_{21}(z)\Big)},
\end{align*} where $c_1=\frac{A_{11}(\hat{z}) - \hat{z}}{A_{11}(\hat{z}) + A_{12}(\hat{z})-\hat{z}}, c_2=\frac{A_{22}(\hat{z}) - \hat{z}}{A_{21}(\hat{z}) + A_{22}(\hat{z})-\hat{z}}$.\\

After simplification, we can write $F(z)$ as
\begin{align}
F(z)=\frac{(1-\rho)(B(z)-1)\Big(c_1 z A_{1}(z)+c_2 z A_{2}(z)+A_{12}(z)A_{21}(z)-A_{11}(z)A_{22}(z) \Big)}{\E[B]\Big((z-A_{11}(z))(z-A_{22}(z))-A_{12}(z)A_{21}(z)\Big)} \label{F_2(z)}.
\end{align}
Let us now determine the expected number of customers $\E[X]=F'(1)$.\\

 After differentiating $F(z)$ w.r.t. $z$ and taking the limit $z \rightarrow 1$, we get
 \begin{align}
\E[X]=&\frac{\rho}{2}+\frac{\text{Var}(A)}{2(1-\rho)}+\frac{\E[B(B-1)]}{2\E[B]} \nonumber\\
&+\frac{-\rho +\E[B](f_1(0)\alpha_1+ f_2(0)\alpha_2)+\rho(\alpha_{11}+\alpha_{22})+\alpha_{12}\alpha_{21}-\alpha_{11}\alpha_{22}}{(P_{12}+P_{21})(1-\rho)}.\label{expX}
\end{align}

For the transient distribution it follows from \eqref{twentysevenAA} that
 \begin{align}
 f_1(r,z)=\frac{z^{X_0+1}\Big(z\P(J_1=1)+r(A_{21}(z)\P(J_1=2)-A_{22}(z)\P(J_1=1))\Big)}{\Big(z-rA_{11}(z)\Big)\Big(z-rA_{22}(z)\Big)-r^2A_{12}(z)A_{21}(z)}\nonumber\\
 +\frac{rz(B(z)-1)\sum_{i=1}^{2}A_{i1}(z)f_i(r,0)}{\Big(z-rA_{11}(z)\Big)\Big(z-rA_{22}(z)\Big)-r^2A_{12}(z)A_{21}(z)}\nonumber\\
 +\frac{r^2(B(z)-1)\Big(A_{12}(z)A_{21}(z)-A_{11}(z)A_{22}(z)\Big)f_1(r,0)}{\Big(z-rA_{11}(z)\Big)\Big(z-rA_{22}(z)\Big)-r^2A_{12}(z)A_{21}(z)},\label{f1rz}
 \end{align}
 \begin{align}
 f_2(r,z)=\frac{z^{X_0+1}\Big(z\P(J_1=2)+r(A_{12}(z)\P(J_1=1)-A_{11}(z)\P(J_1=2))\Big)}{\Big(z-rA_{11}(z)\Big)\Big(z-rA_{22}(z)\Big)-r^2A_{12}(z)A_{21}(z)}\nonumber\\
 +\frac{rz(B(z)-1)\sum_{i=1}^{2}A_{i2}(z)f_i(r,0)}{\Big(z-rA_{11}(z)\Big)\Big(z-rA_{22}(z)\Big)-r^2A_{12}(z)A_{21}(z)}\nonumber\\
 +\frac{r^2(B(z)-1)\Big(A_{12}(z)A_{21}(z)-A_{11}(z)A_{22}(z)\Big)f_2(r,0)}{\Big(z-rA_{11}(z)\Big)\Big(z-rA_{22}(z)\Big)-r^2A_{12}(z)A_{21}(z)},\label{f2rz}
 \end{align}    where
 \begin{align}
f_1(r,0)=\frac{\left(-\hat{z}_1^{X_0}(\hat{B}^{(2)}-1)\hat{A}_{21}^{(2 )}(\hat{z}_1-r\hat{A}_{22}^{(1 )})+\hat{z}_2^{X_0}(\hat{B}^{(1)}-1)\hat{A}_{21}^{(1 )}(\hat{z}_2-r\hat{A}_{22}^{(2 )}) \right)\P(J_1=1)}{(\hat{B}^{(1)}-1)(\hat{B}^{(2)}-1)\left(\hat{A}_{21}^{(2 )}(\hat{z}_1-r\hat{A}_{22}^{(1 )})-\hat{A}_{21}^{(1 )}(\hat{z}_2-r\hat{A}_{22}^{(2 )})\right)}\nonumber\\
+ \frac{r\left(\hat{z}_2^{X_0}(\hat{B}^{(1)}-1)-\hat{z}_1^{X_0}(\hat{B}^{(2)}-1)\right)\hat{A}_{21}^{(1 )}\hat{A}_{21}^{(2 )}\P(J_1=2)}{(\hat{B}^{(1)}-1)(\hat{B}^{(2)}-1)\left(\hat{A}_{21}^{(2 )}(\hat{z}_1-r\hat{A}_{22}^{(1 )})-\hat{A}_{21}^{(1 )}(\hat{z}_2-r\hat{A}_{22}^{(2 )})\right)} ,\label{firstvalue}
 \end{align}
\begin{align}
 f_2(r,0)=\frac{1}{r}\frac{\left(\hat{z}_1^{X_0} (\hat{B}^{(2)}-1)-\hat{z}_2^{X_0}(\hat{B}^{(1)}-1)\right)\left(\hat{z}_1-r\hat{A}_{22}^{(1 )}\right)\left(\hat{z}_2-r\hat{A}_{22}^{(2 )}\right)\P(J_1=1)}{(\hat{B}^{(1)}-1)(\hat{B}^{(2)}-1)\left(\hat{A}_{21}^{(2 )}(\hat{z}_1-r\hat{A}_{22}^{(1 )})-\hat{A}_{21}^{(1 )}(\hat{z}_2-r\hat{A}_{22}^{(2 )})\right)}\nonumber\\
+ \frac{\left(-\hat{z}_2^{X_0}(\hat{B}^{(1)}-1)\hat{A}_{21}^{(2 )}(\hat{z}_1-r\hat{A}_{22}^{(1 )})+\hat{z}_1^{X_0}(\hat{B}^{(2)}-1)\hat{A}_{21}^{(1 )}(\hat{z}_2-r\hat{A}_{22}^{(2 )}) \right)\P(J_1=2)}{(\hat{B}^{(1)}-1)(\hat{B}^{(2)}-1)\left(\hat{A}_{21}^{(2 )}(\hat{z}_1-r\hat{A}_{22}^{(1 )})-\hat{A}_{21}^{(1 )}(\hat{z}_2-r\hat{A}_{22}^{(2 )})\right)}, \label{secondvalue}
 \end{align} $z=\hat{z}_1$ and $z=\hat{z}_2$ are the zeros in the unit disc $|z|<1$ of $\Big(z-rA_{11}(z)\Big)\Big(z-rA_{22}(z)\Big)-r^2A_{12}(z)A_{21}(z)$ and $\hat{A}_{ij}^{(1 )}:=A_{ij}(\hat{z}_1)$, $\hat{A}_{ij}^{(2 )}:=A_{ij}(\hat{z}_2),\hat{B}^{(i)}:=B(\hat{z}_i)$ for $i,j=1,2$.\\

\begin{remark}\label{remark:1}
It can be observed that the first three terms in the right-hand-side of Equation \eqref{expX} are exactly equal to the mean queue length at departure epochs of the standard $M^{X}/G/1$ queue without dependencies, cf.~Gaver~\cite{gaver59} and Cohen~\cite[Section~III.2.3]{cohen}, and the remaining term appears due to the dependent service times.

\end{remark}
\begin{remark}\label{remark:2} It can be shown, after some straightforward but tedious algebraic manipulations, that the queue-length distribution in the system considered in the present paper also reduces to the distribution of the number of customers in an $M^X/G/1$ queuing model if $A_1(z)=A_2(z)=A(z)$, again cf.~Gaver~\cite{gaver59} and Cohen~\cite[Section~III.2.3]{cohen}. Similarly, we can also prove that the \emph{expected} number of customers in the system considered in the present paper is equal to the \emph{expected} number of customers in the corresponding $M^X/G/1$ queuing model  if $\alpha_1=\alpha_2=\E[A]$.\\
\end{remark}

\section{Numerical results}
\label{sec:num}
In this section, we present four numerical examples in order to get more insight in the consequences of introducing dependencies between the service times of consecutive customers. For simplicity, we restrict ourselves to two customer types ($N=2$).  In all four examples we assume that the overall batch arrival process is a Poisson process with rate $\lambda$ and the load $\rho$ equals $\frac34$.

\subsection{Example 1}

In this example we consider an almost symmetric system, with $\P(J=1)=\P(J=2)=\frac{1}{2}$ and  $\alpha_{ij}=\frac{3}{8}$ for $\forall i,j = 1,2$. It follows that $\E[A]=\frac{3}{4}$, $P_{11}=P_{22}$ and we shall vary $P_{11}$. The batch sizes are geometrically distributed with
\[
\P(B=k)=p^{k-1}(1-p), \qquad k=1,2,\dots
\]
We take $p=3/4$, resulting in a mean batch size of $\E[B]=4$.
The conditional service times are respectively exponential and Erlang distributed random variables, with
\begin{align*}
G_{ij}(x) &=  \left(1-\sum_{m=0}^{k_j-1} \frac{(\mu_{ij} x)^m}{m!}e^{-\mu_{ij}x}\right)P_{ij},
\end{align*}
for $\mu_{ij} > 0,$ $i,j=1, 2$. In this example we will take an Erlang distribution with four phases. If we define
\[
k_j=\begin{cases}
1 & \qquad \text{if }j=1,\\
4 & \qquad \text{if }j=2,
\end{cases}
\]
we can use Equation \eqref{GijzA} to obtain
\[
A_{ij}(z)=P_{ij}\left(\frac{\mu_{ij}}{\lambda(1-B(z))+\mu_{ij}}\right)^{k_j},
\]
for $i=1,2$ and $j=1,2$.
\\

The variance of the number of arrivals during one arbitrary service time, written as a function of $P_{11}$, directly follows. For  $0<P_{11} <1$,
\[
\text{Var}(A) = \frac{75}{16}
+\frac{117}{512 (1-P_{11}) P_{11}}.
\]
We observe that $\alpha_1=\alpha_2$, but $A_1(z)\neq A_2(z)$. From Remark \ref{remark:2}, we know that the mean queue length in our model is equal to the mean queue length of a standard $M^X/G/1$ queue, 
but for higher moments of the queue length, this equality is not true unless we can construct a case with $A_1(z)=A_2(z)$. This is confirmed by Table \ref{tbl:exmg1}, which depicts numerical values for the means and variances of the queue lengths in our model and in the corresponding $M^X/G/1$ queue. Indeed, the mean queue lengths of both systems are equal, whereas the variances of the queue lengths are only equal in the case $P_{11}=\frac12$, where $A_1(z)=A_2(z)$.
Since $\alpha_1=\alpha_2$, we immediately conclude that the mean queue length and the variance of $A$ are minimal when $P_{11}=1/2$ (see Remark~\ref{remark:2}).
\begin{table}[ht!]
\begin{center}
$$
\begin{array}{|c|c|cc|}
\hline
P_{11} & \E[X]= \E[X^{M^X/G/1}] & \text{Var}(X) & \text{Var}(X^{M^X/G/1}) \\
\hline
0.1 &    17.8281 & 374.4642 & 374.4631 \\
0.3 &    14.9263 & 237.6202 & 237.6198 \\
0.5 &    14.5781 & 223.8303 & 223.8303 \\
0.7 &    14.9263 & 237.6184 & 237.6198 \\
0.9 &    17.8281 & 374.4185 & 374.4631 \\
\hline
\end{array}
$$
\end{center}
\caption{Means and variances of $X$ and $X^{M^X/G/1}$ for various values of $P_{11}$ in Example 1.}
\label{tbl:exmg1}
\end{table}

\subsection{Example 2}

In this example we take a similar setting as in the previous example, but we make two adjustments. First, for even more simplicity, we assume that all conditional service times are exponentially distributed, i.e.,
\[
G_{ij}(x)=(1-e^{-\mu_{ij}x})P_{ij}, \qquad \forall i,j=1,2.
\]
Secondly, we take $\alpha_{11}=\alpha_{12}=\frac{1}{2}$ and $\alpha_{21}=\alpha_{22}=\frac{1}{4}$. As in the previous example, we let $\P(J=1)=\P(J=2)=\frac{1}{2}$. We observe that the difference with Example 1 is that all conditional service-time distributions are exponential now, but with different parameters. Moreover, in this model $\alpha_1 \neq \alpha_2$.

An interesting question is, how the mean queue length and the variance of the number of arrivals during an arbitrary service time are related. Since $\alpha_1 \neq \alpha_2$, the setting of Remark \ref{remark:2} does not apply. In Figure \ref{fig:EXvarAexample2} we show $\E[X]$ and $\text{Var}(A)$ plotted versus $P_{11}$. When studying the two plots carefully, one can see that the plots are not completely symmetric, which is obviously caused by the asymmetric service times. However, another observation that is not visible to the human eye, is that the minima of both plots are \emph{not} attained at the same value of $P_{11}$. It can be shown analytically, that the variance of $A$ is minimal at exactly $P_{11}=1/2$, and numerically, that $\E[X]$ is minimal for $P_{11}\approx 0.500411$. Although this is a small difference, it means that this system exhibits an interesting, rare feature: it is possible to obtain a \emph{smaller} mean queue length by having a \emph{greater} variance in the number of arrivals during one service time. In Example 3 we will create a setting in which this effect is even bigger.

\begin{figure}[ht!]
\parbox{0.45\textwidth}{\centering
\includegraphics[width=\linewidth]{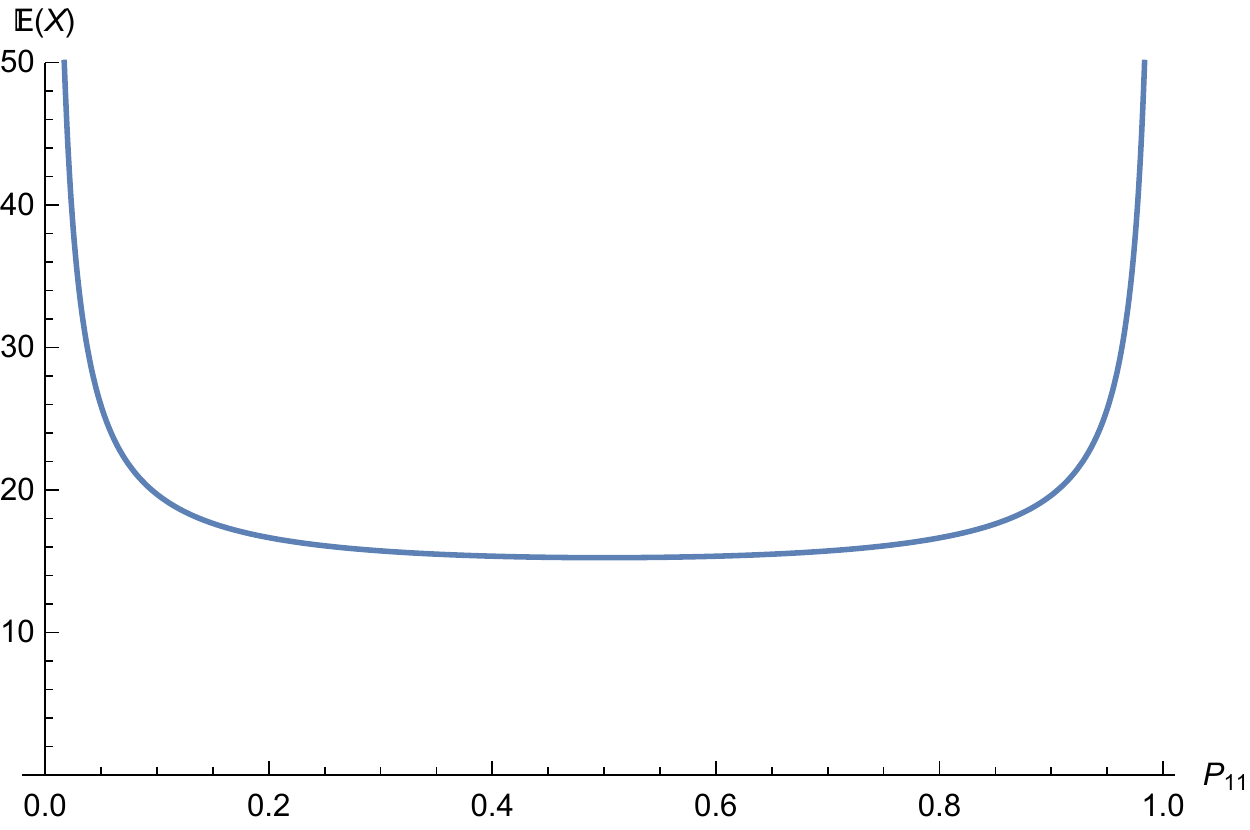}\\(a) Mean queue length
}\hfill
\parbox{0.45\textwidth}{\centering
\includegraphics[width=\linewidth]{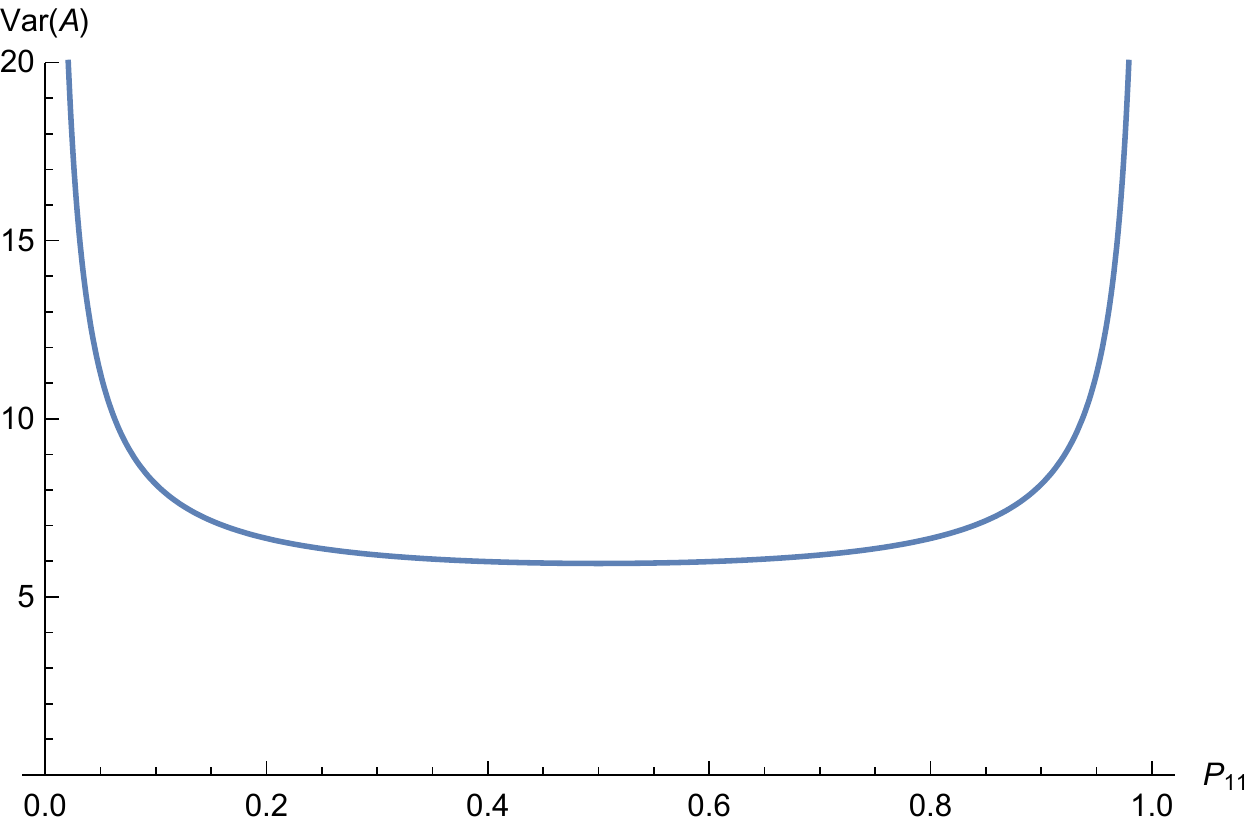}\\(b) Variance of $A$
}
\caption{The mean queue length $\E[X]$ and the variance of $A$ in Example 2.}
\label{fig:EXvarAexample2}
\end{figure}



From Figures \ref{fig:EXvarAexample2}(a) and (b), we can observe that, except for the small region where $0.5<P_{11}< 0.500411$, the expected number of customers is increasing when the variance of the number of arrivals during a customer service time is increasing and conversely. This means that bigger variance of the number of arrivals implies a larger expected number of customers. This also implies that the expected number of customers can grow beyond any bound in a stable system due to the very large variance of the number of arrivals during one service time. This scenario occurs when $P_{11}$ tends to $0$ or $1$ in Figure \ref{fig:EXvarAexample2}. Therefore, we can observe dependencies when  $P_{11}$ or $(1-P_{11})$ is small. Otherwise, $\E[X]$ and $\text{Var}(A)$ appear to be rather insensitive to the value of $P_{11}$.

Of course, the reason for the large variance in the number of arrivals during a customer service time lies in the dependence. When, e.g., $P_{11}=P_{22}$ is very small, services alternate for a long time between $\text{exp}(\mu_{12})$ and $\text{exp}(\mu_{21})$ services with small mean; rarely is there an $\text{exp}(\mu_{11})$  or $\text{exp}(\mu_{22})$  service which has a huge mean.

\subsection{Example 3}
 Once again, we assume that the conditional service times are exponentially distributed, but in this example we choose less symmetric settings. Let $\P(J=1)=\frac{7}{16}, \P(J=2)=\frac{9}{16}, \alpha_{11}=\alpha_{12}=\alpha_{21}=\frac{3}{20}$ and $ \alpha_{22}=\frac{19}{20}$. From these settings we obtain $P_{21}=\frac{7}{9}P_{12}$, $\alpha_1=0.3$, and $\alpha_2=1.1$. The interesting phenomenon observed in Example 2, is also taking place here. In fact, in this example there is a bigger difference between the value of $P_{11}$ for which the mean queue length is minimal ($P_{11} \approx 0.65$), and the value resulting in a minimum variance of the number of arrivals during an arbitrary service time ($P_{11} \approx 0.788$). More details can be found in Table~\ref{tbl:exvara}. The interesting region is obviously $0.650 < P_{11} < 0.788$, because in this region we know that an increase in $\text{Var}(A)$ results in a decrease in $\E[X]$. This is illustrated even better in Figure \ref{fig:exvara}, where $\text{Var}(A)$ and $\E[X]$ are plotted against each other, for varying values of $P_{11}$.

\begin{table}[ht!]
\begin{center}
$$
\begin{array}{|c|cc|}
\hline
P_{11} & \E[X] & \text{Var}(A) \\
\hline
 0.100 & 20.377 & 8.327 \\
 0.300 & 17.931 & 7.056 \\
 0.500 & 16.969 & 6.493 \\
\mathbf{ 0.650} & \mathbf{16.747} & 6.263 \\
 0.700 & 16.780 & 6.214 \\
 \mathbf{0.788} & 17.060 & \mathbf{6.175} \\
 0.900 & 18.587 & 6.333 \\
\hline
\end{array}
$$
\end{center}
\caption{Mean queue length and variance of the number of arrivals during an arbitrary service time, for various values of $P_{11}$ in Example 3.}
\label{tbl:exvara}
\end{table}

\begin{figure}[!ht]
\parbox{0.45\linewidth}{
\centering
\includegraphics[width=\linewidth]{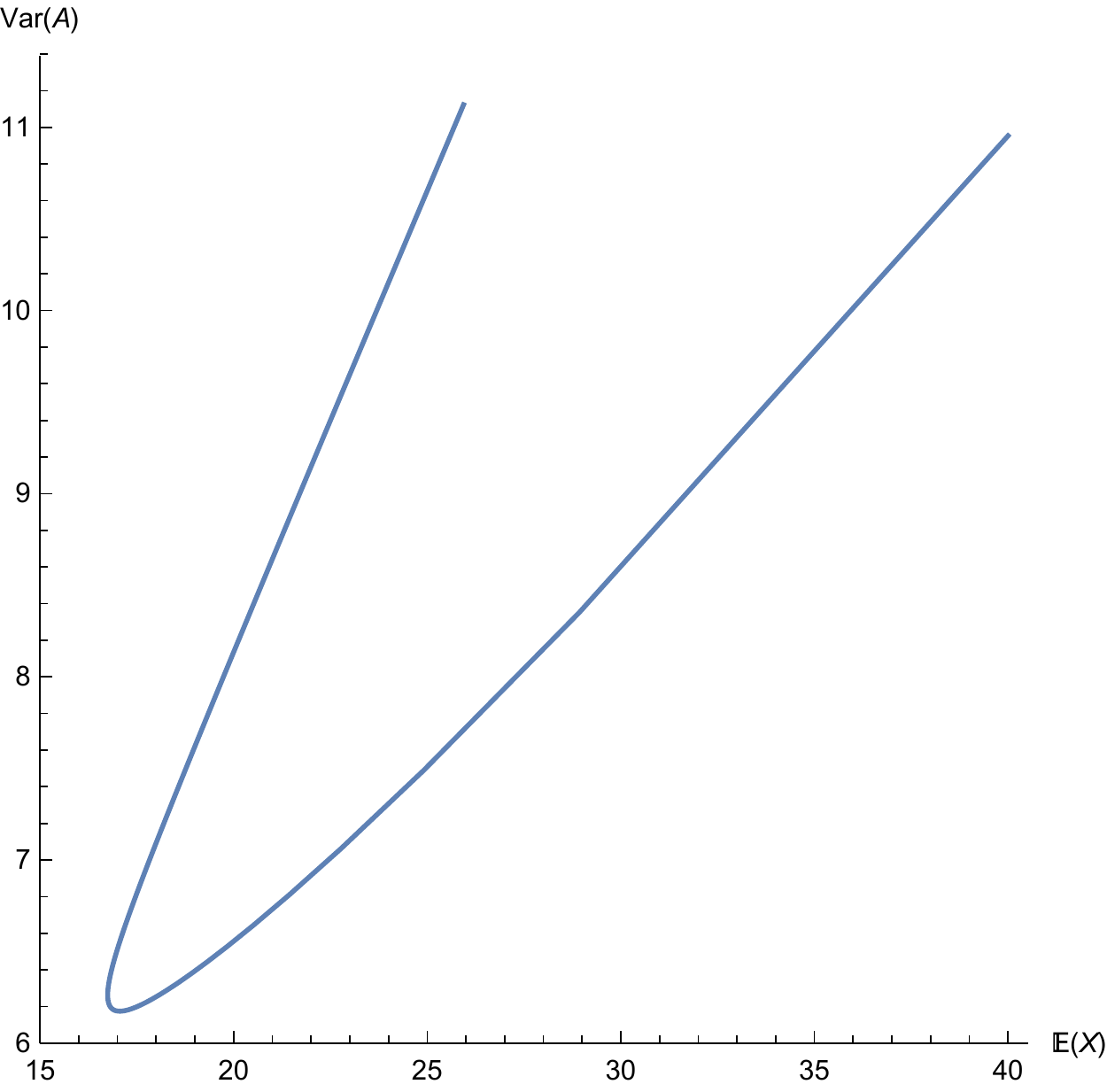}\\
(a) $0.01 < P_{11} < 0.99$
}\hfill
\parbox{0.45\linewidth}{
\centering
\includegraphics[width=\linewidth]{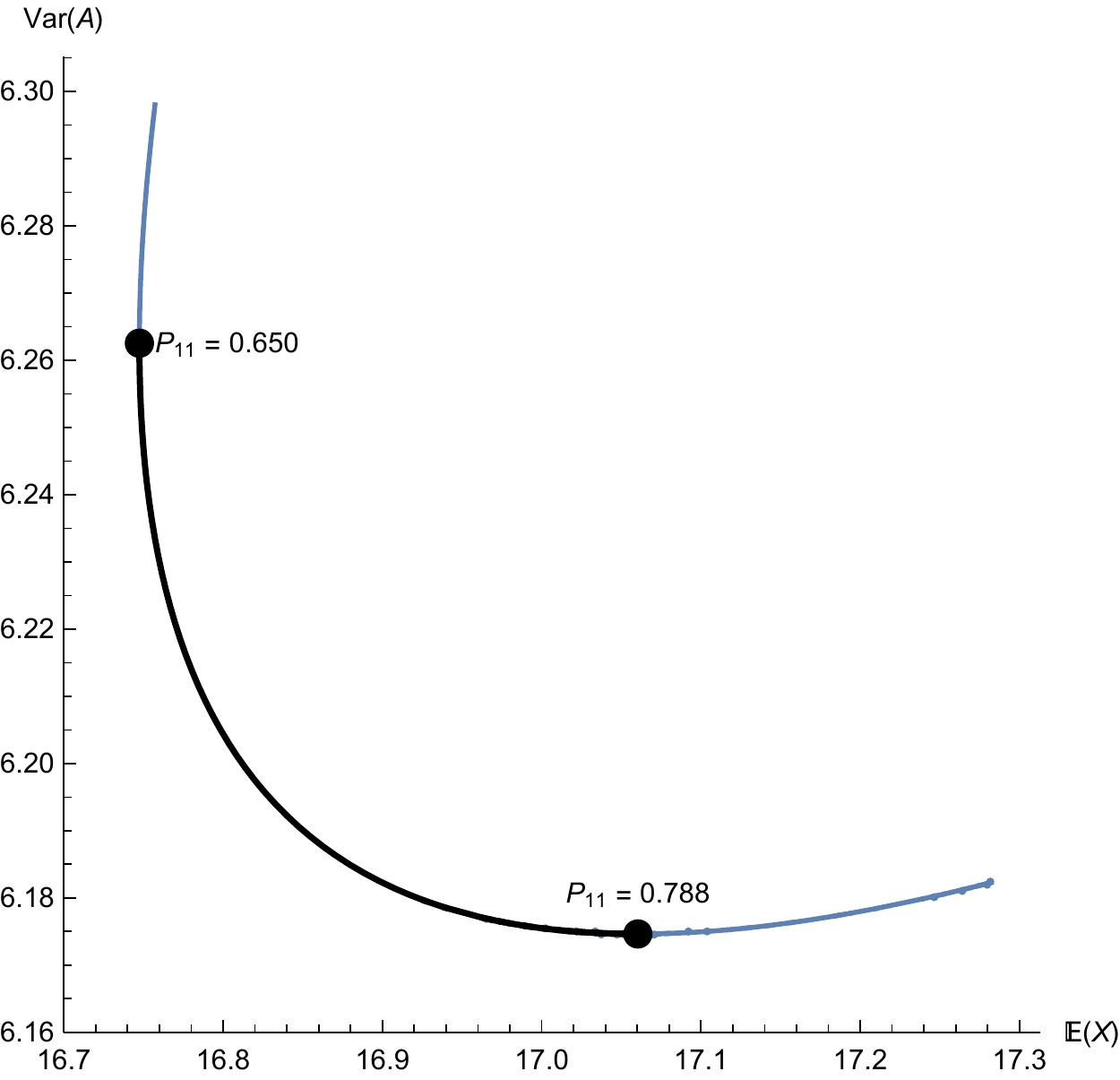}\\
 (b) $0.62 < P_{11} < 0.82$
}
\caption{The variance of the number of arrivals versus the expected number of customers during an arbitrary customer service time. This implicit plot is obtained by varying $P_{11}$. Figure (b) is a zoomed in version of Figure (a).}
\label{fig:exvara}
\end{figure}

\subsection{Example 4: Transient-state analysis} We return to the system in Example 2, but now we study the transient
analysis. In
this example we start with an empty system, $\E[z^{X_{0}}]=1$, and set
$P_{11}=1/10$. Next, we repeatedly apply Equation \eqref{twentysixAA} to
express $\E[z^{X_{n}}]$ in terms of $\E[z^{X_{n-1}}]$.
We have taken four different distributions for the conditional service
times, namely exponential, gamma with shape parameter $1/2$, gamma with
shape parameter 5, and deterministic.
The results are shown in Figure \ref{fig:ex4transient}, where we depict
the mean queue length after the departure of the $n$th customer, for
$n=0,1,2,\dots,200$.
In this example, it can clearly be seen that service-time distributions
with higher coefficients of variation result in longer queues. Also, it
seems to take longer to reach steady state. For completeness, we give
the steady-state mean queue lengths for the four systems below:
\begin{center}
\begin{tabular}{|c|c|c|c|c|}
\hline
Distribution & Deterministic & Gamma $5$  & Exponential & Gamma $1/2$\\
\hline
$\E[X]$ & 16.224 & 16.918 & 19.696 & 23.168   \\ 
\hline
\end{tabular}
\end{center}

\begin{figure}[!ht]
\centering
\includegraphics[width=0.6\linewidth]{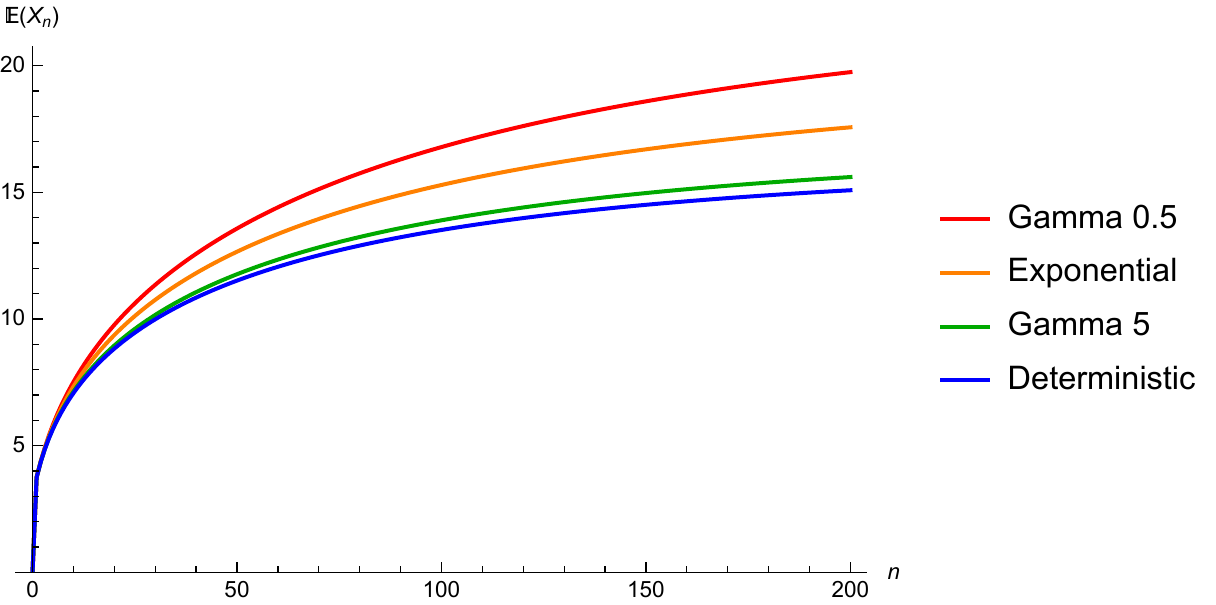}
\caption{Numerical example 4: Transient mean queue-length analysis.}
\label{fig:ex4transient}
\end{figure}

\section*{Acknowledgments} The research of Abhishek and Onno Boxma  is partly funded by NWO Gravitation project {\sc Networks}, grant number 024.002.003. The authors thank Michel Mandjes (University of Amsterdam) for helpful discussions.

%

\end{document}